\documentclass{conm-p-l}
\usepackage{amsfonts,amssymb,amscd,amsmath,enumerate}

\copyrightinfo{2000}{American Mathematical Society}

\newcommand{\qnd}{\quad\text{and}\quad}

\newcommand{\ssstyle}{\scriptscriptstyle}

\newcommand{\ges}{\operatorname{\ssstyle\geqslant}}
\newcommand{\les}{\operatorname{\ssstyle\leqslant}}

\newcommand{\sgn}[1]{(-1)^{|#1|}}

\newcommand{\col}{\colon}
\newcommand{\dd}{\partial}

\newcommand{\wh}{\widehat}

\newcommand{\arto}{\xrightarrow}

\newcommand{\Ker}{\operatorname{Ker}}

\newcommand{\ass}{\operatorname{Ass}}

\newcommand{\depth}[2]{\operatorname{depth}\left(#1,#2\right)}
\newcommand{\mdepth}{\operatorname{depth}}
\newcommand{\Rdepth}[3]{{\operatorname{depth}_{#1}(#2,#3)}}

\newcommand{\fdim}{\operatorname{fd}}

\newcommand{\pdim}{\operatorname{pd}}

\newcommand{\rad}{\operatorname{rad}}
\newcommand{\spec}{\operatorname{Spec}}

\newcommand{\mwidth}{\operatorname{width}}
\newcommand{\Rwidth}[3]{{\operatorname{width}_{#1}(#2,#3)}}
\newcommand{\V}[1]{\operatorname{V}(#1)}

\newcommand{\amp}{\operatorname{amp}}

\newcommand{\ED}[3]{{}^{#1}\!{d}_{#2,#3}}

\newcommand{\EH}[3]{{}^{#1}\!\operatorname{E}_{#2,#3}}
\newcommand{\Ext}[4]{\operatorname{Ext}^{#1}_{#2}\left(#3,#4\right){}}

\newcommand{\hh}[1]{\operatorname{H}(#1)}
\newcommand{\HH}[2]{\operatorname{H}_{#1}(#2)}

\newcommand{\Hom}[3]{\operatorname{Hom}_{#1}\left(#2,#3\right){}}
\newcommand{\Rhom}[3]{\operatorname{{\bf R}Hom}_{#1}\left(#2,#3\right){}}
\newcommand{\shift}[2]{{\operatorname{\scriptstyle\Sigma}^{#1}\!{#2}}}
\newcommand{\ltrunc}[2]{\tau_{\les #1}\left(#2\right)}
\newcommand{\rtrunc}[2]{\tau_{\ges #1}\left(#2\right)}
\newcommand{\Tor}[4]{\operatorname{Tor}_{#1}^{#2}\left(#3,#4\right){}}
\newcommand{\Ltor}[3]{{#1}\otimes^{\bf L}_{#2}{#3}}

\newcommand{\frob}[2]{({}^{\vf^{#1}}\!{#2})}

\newcommand{\dcat}[1]{{\mathcal D}(#1)}
\newcommand{\lho}[3]{\operatorname{H}_{#1}^{#2}\left(#3\right)}
\newcommand{\Lho}[2]{{\bf L}{\Lambda}_{#1}\left(#2\right)}
\newcommand{\lch}[3]{\operatorname{H}^{#1}_{#2}\left(#3\right)}
\newcommand{\Lch}[2]{{\bf R}{\Gamma}_{#1}\left(#2\right)}

\newcommand{\BZ}{{\mathbb Z}}

\newcommand{\bsx}{{\boldsymbol x}}

\newcommand{\5}{}
\newcommand{\fa}{{\mathfrak a}}
\newcommand{\fb}{{\mathfrak b}}

\newcommand{\fm}{{\mathfrak m}}

\newcommand{\fp}{{\mathfrak p}}
\newcommand{\fq}{{\mathfrak q}}

\newcommand{\vf}{{\varphi}}

\swapnumbers

\theoremstyle{plain}
\newtheorem{theorem}{Theorem}[section]

\newtheorem*{one}{Theorem I}
\newtheorem*{two}{Theorem II}
\newtheorem*{three}{Theorem III}
\newtheorem*{four}{Theorem IV}
\newtheorem*{five}{Theorem V}

\newtheorem{proposition}[theorem]{Proposition}

\newtheorem*{Lemma}{Lemma}

\newtheorem{noname}[theorem]{}

\theoremstyle{definition}

\newtheorem{definition}[theorem]{Definition}

\newtheorem{example}[theorem]{Example}

\newtheorem{chunk}[theorem]{}

\theoremstyle{remark}

\newtheorem*{Remark}{Remark}

\newtheorem*{caveat}{Caveat}

\begin{document}

\title[Depth and amplitude]{Depth and amplitude for unbounded complexes}

\author{Hans--Bj\o rn Foxby}
\address{Matematisk Afdeling, K\o benhavns Universitet,
Universitetsparken 5,
\endgraf DK--2100 K\o benhavn \O, Denmark}
\email{foxby@math.ku.dk}

\author{Srikanth Iyengar}
\address{Mathematics Department, University of Missouri, Columbia, MO
65211, USA} \email{iyengar@math.missouri.edu}

\thanks{H.--B.~F.~was partially supported by the Danish Natural Science
Research Council.\endgraf S.~I.~was supported by a grant from the E.P.S.R.C}


\begin{abstract}
  We prove that over a commutative noetherian ring the three approaches to introducing
  depth for complexes: via Koszul homology, via Ext modules, and via local cohomology, all
  yield the same invariant. Using this result, we establish a far reaching generalization
  of the classical Auslander-Buchsbaum formula for the depth of finitely generated modules
  of finite projective dimension. We extend also Iversen's amplitude inequality to
  unbounded complexes.  As a corollary we deduce: Given a local homomorphism $Q\to R$, if
  there is a non-zero finitely generated $R$-module that has finite flat dimension both over
  $Q$ and over $R$, then the flat dimension of $R$ over $Q$ is finite. This last result
  yields a module theoretic extension of a characterization of regular local rings in
  characteristic $p$ due to Kunz and Rodicio
\end{abstract}

\maketitle

\section*{Introduction}
This paper concerns the theory of depth for complexes over commutative noetherian
rings. Let ${\fa}$ be an ideal in such a ring $R$ and $K$ the Koszul complex on a
finite set of $n$ generators for ${\fa}$. It is well known that the following numbers
are equal when $M$ is a finitely generated $R$-module:
\begin{itemize}
\item $n - \sup\{\,\ell\in\BZ\mid\HH{\ell}{K\otimes_RM}\ne0\,\}$;
\item $\inf\{\,\ell\in\BZ\mid\Ext\ell R{R/{\fa}}M\ne0\,\}$;
\item $\inf\{\,\ell\in\BZ\mid\lch\ell {\fa}M\ne0\,\}$, where
$\lch*{\fa}M$ is the ${\fa}$-local cohomology of $M$.
\end{itemize}

If $R$ is local, then this common value coincides with the length of the longest
$M$-regular sequence in ${\fa}$, that is to say, the ${\fa}$-depth of $M$. Since each
of the quantities displayed above is meaningful even when $M$ is a {\em complex} of
$R$-modules, they can, and have been, taken as a starting point for building a theory
of depth for complexes. These have proved to be of immense utility even in dealing with
problems concerning modules alone.

The natural question arises: {\em Do these three approaches yield the same invariant?}
Results of Foxby and Iyengar settle this question in the affirmative for about a third of
the complexes.  Namely, those whose homology is bounded above, noteworthy among these
being complexes with bounded homology and that includes also all modules, finitely
generated or not.

One of our main results answers the question with an unqualified: Yes.

\begin{one}
For any complex of $R$-modules $M$ the numbers obtained from the three
formulas above coincide.

\end{one}

Thus one may now speak of {\em the} depth of a complex without having to specify which
of the preceding formulas is being used to define it. This flexibility in computing
depth allows one to extend many  familiar results on depth for modules to identical
statements for arbitrary complexes.

For instance, calculating depth via the vanishing of Ext leads to the following theorem
concerning the depth (with respect to the maximal ideal) of complexes over local rings.
It is a vast generalization of the classical Auslander-Buchsbaum equality: $\mdepth
R=\mdepth_RP + \pdim_RP$ for any finitely generated $R$-module $P$ whose projective
dimension $\pdim_RP$ is finite.

\begin{two}
Let $R$ be a local ring and let $M$ and $P$ be complexes of $R$-modules. If
$\,\pdim_RP$ is finite and \,$\hh P$ is nonzero and finitely generated, then
\[
\mdepth_RM = \mdepth_R(\Ltor MRP) + \pdim_RP\,.
\]
 \end{two}
Here $\Ltor MRP$ stands for $M\otimes_RF$, where $F$ is any finite free resolution of
$P$.

Our proof of the next result makes critical use of the characterization of depth via
Koszul homology and also of the preceding theorem. Recall that the {\em
amplitude} of $M$ is the number 
\[
\amp M = \sup\{\ell\mid \HH\ell M\ne0\}-\inf\{\ell\mid \HH\ell M\ne0\}\,.
\]
We write $\fdim_RF$ for the flat dimension of $F$ over $R$.

\begin{three}
Let $R\to S$ be a local homomorphism and $F$ a complex of $S$-modules with $\hh F$
nonzero and finitely generated. If $\,\fdim_RF$ is finite, then for any complex of
$R$-modules $M$ with $\hh M$ degreewise finite one has
\[
\amp M \leq \amp\left(\Ltor MRF\right)\,.
\]
 \end{three}
 
 The first step in our proof of this theorem is a reduction to the case where $S=R$; this
 is readily accomplished by using the theory of Cohen factorizations of Avramov, Foxby,
 and Herzog \cite{AFH}.  Now, when $\hh M$ is bounded the result we seek is precisely the
 amplitude inequality, proved by Iversen \cite{Iv1}; it is an avatar of Paul Roberts'
 remarkable New Intersection Theorem \cite{Pr}.  Thus, the new information provided by our
 theorem concerns {\em unbounded} complexes; the issues that come into play in proving it
 are of a different nature and not as involved. Nevertheless, as the following corollary
 demonstrates, it too has its uses.

\begin{four}
Let $Q\to R\to S$ be local homomorphisms and $F$ a nonzero finitely generated
$S$-module. Assume that $\fdim_RF<\infty$.  Then
\[
\fdim_QR\leq \fdim_QF\leq \fdim_QR + \fdim_RF\,.
\]
In particular, $\fdim_QR$ and $\fdim_QF$ are finite simultaneously.
 \end{four}

The inequality on the right is classical. Since regularity descends along maps of
finite flat dimension, the one on the left implies that if $R$ is regular and
$\fdim_QF$ is finite, then $Q$ is regular. The special case where $R=S$ of this latter
result had been discovered by Apassov \cite{Ap}.

Theorem IV yields also the following characterization of regular rings of characteristic
$p$ that extends those of Kunz \cite{Ku} and Rodicio \cite{Ro}. In its statement,
${}^{\vf^n}\!F$ denotes the $R$ module structure on $F$ induced by restriction of scalars
along $\vf^n$.

\begin{five}
Let $R$ be a local ring of characteristic $p$, and let $\vf\col R\to R$ be the
Frobenius endomorphism of $R$. The following conditions are equivalent.
\begin{enumerate}[{\quad\rm(a)}]
\item $R$ is regular; 
\item $\vf^n$ is flat for each integer $n\geq 0$;
\item there exists a positive integer $n$ and a nonzero finitely generated $R$-module $F$
  such that both $\fdim_RF$ and $\fdim_R\frob nF$ are finite.
\end{enumerate}
\end{five}

The results in this paper have further applications to the study of the Frobenius
endomorphism. These will be dealt with in \cite{AIM} and \cite{IW}.

Theorem IV has implications also for flat dimensions of homomorphisms.  In the sequel, the
{\em flat dimension} of any homomorphism $\psi\col R\to S$ is the flat dimension of $S$ over
$R$; it is denoted $\fdim\psi$.  

Consider a diagram of local homomorphisms $Q\arto{\vf}R\arto{\psi}S$.  As is well known,
if both $\vf$ and $\psi$ have finite flat dimension, then so does the composition
$\psi\circ\vf$. The `factorization problem' for maps of finite flat dimension reads: 
{\em  Suppose that $\fdim (\psi\circ\vf)$ is finite. What can be said about $\fdim\vf$ and
$\fdim\psi$?}  We are now in a position to give a complete answer to this question.

Indeed, if $\fdim \psi$ is finite, then $\fdim\vf$ is also finite; this is by Theorem
IV.

On the other hand, it could happen that $\fdim\vf$ is finite, but $\fdim\psi$ is not:
Let $(Q,\fm,k)$ be a regular local ring and $f$ a nonzero element in $\fm^2$. Then the
diagram $Q\to Q/(f)\to k$ furnishes the requisite counterexample.

Theorems I and II are proved in Section \ref{depth}. That section establishes also the
basic properties of depth for complexes. Most of these generalize results of Iversen,
Foxby, and Iyengar.  The proofs of Theorems III, IV, and V are contained in Section
\ref{amplitude}. Section \ref{width} concerns certain results on width for complexes
analogous to those on depth described in Section \ref{depth}.

Most of our calculations occur in the derived category of modules for they are most
transparent there. For this reason, in Section \ref{basis}, we collect a few useful
facts concerning the same, as well as results on local cohomology and derived
completions required in this work.

We stumbled upon the main results described in this paper at the conference in
Grenoble.\footnote{Between the morning and afternoon sessions. Not during.} It is a
pleasure to thank the organizers for inviting us, and also for giving us an opportunity
to present our results in this forum.

\section{Homological algebra}\label{basis}
Let $R$ be a ring; the standing assumption in this paper is that rings are commutative
and noetherian. We are concerned with complexes of the form
\[
\cdots \to M_{\ell+1}\arto{\dd_{\ell+1}} M_\ell \arto{\dd_\ell} M_{\ell-1}\to\cdots
\]
Given a complex $M$ of $R$-modules set
\[
\sup M = \sup\{\ell\mid \HH \ell M \ne  0\} \qnd
  \inf M = \inf\{\ell\mid \HH \ell M\ne 0\}\,,
\]
with the understanding that $\sup\emptyset=-\infty$ and $\inf\emptyset=\infty$. When
$\sup M<\infty$ (respectively, $\inf M>\infty$), one says that $\hh M$ is {\em bounded
above} (respectively, {\em bounded below}). Naturally, $\hh M$ is {\em bounded} if it
is bounded both above and below.  We label $\hh M$ {\em degreewise finite} if the
$R$-module $\HH iM$ is finitely generated for each integer $i$; it is labelled {\em
finite} if, in addition, it is bounded.

Associated to any complex $M$ and integers $i,s$ are the {\em truncated} complexes
\begin{align*}
&\rtrunc iM \col\quad \cdots \to M_{i+2}\to M_{i+1} \to \Ker(\dd_i) \to
0\,,\\
&\ltrunc sM \col\quad 0\to \frac{M_s}{\dd(M_{s+1})}\to M_{s-1}\to M_{s-2}\to \cdots\,.
\end{align*}
The canonical inclusion $\rtrunc iM\to M$ induces an isomorphism in homology in
degrees $\geq i$, whilst the map in homology induced by the surjection $M\to\ltrunc sM$
is an isomorphism in degrees $\leq s$. In particular, if $i\leq \inf M$ (respectively,
$s\geq \sup M$), then the morphism $\rtrunc iM\to M$ (respectively, $M\to \ltrunc sM$)
is a {\em quasiisomorphism}, that is to say, a homology isomorphism.

The {\em $n$th suspension} of $M$ is the complex $\shift nM$ with $(\shift
nM)_i=M_{i-n}$ and differential given by $\dd(m)=\sgn n\dd_{i-n}(m)$ for $m\in(\shift
nM)_i$.

\subsection*{Koszul homology}
For an element $x$ in the ring $R$ denote $K(x)$ the complex $0\to R\arto{x} R \to 0$
concentrated in degrees $0$ and $1$. The {\em Koszul complex} on a sequence
$\bsx=\{\,x_1,\dots,x_n\,\}$ is the complex
\[
K(\bsx) = K(x_1)\otimes_R\cdots \otimes_RK(x_n)\,.
\]

In the sequel we take recourse to the following properties of Koszul complexes.

\begin{noname}\label{les:koszul}
Let $K$ be the Koszul complex on an element $x$ in $R$ and let $M$ be a complex.
There is a long exact sequence
\[
\cdots\to\HH iM \arto{x}\HH iM\to\HH i{K\otimes_RM}\to\HH{i-1}M \arto{x}\cdots\,.
 \]
\end{noname}
Indeed, the canonical inclusion $R=K_0\hookrightarrow K$ extends to a short exact
sequence of complexes of free modules $0\to R\to K \to \shift{}R\to0$. When tensored
with $M$, this induces the exact sequence of complexes
\[
0\to M \to K\otimes_RM \to \shift{}M \to 0\,.
\]
The corresponding homology long exact sequence is the one we seek.

\begin{noname}\label{annihilation:koszul}
Let $K$ be the Koszul complex on a set of elements $\{\,x_1,\dots,x_n\,\}$ in $R$,
and set ${\fa}=(x_1,\dots,x_n)$. For any complex of $R$-modules $M$ one has
\[
{\fa}\5 \hh{K\otimes_RM}=0 \qnd {\fa}\5 \hh{\Hom RKM}=0\,.
\]
\end{noname}
This follows from the fact that multiplication by $x_i$ is homotopic to $0$ on
$K(x_i)$, and so on $K\otimes_RM$ and $\Hom RKM$. Thus, $x_i$ annihilates both
$\hh{K\otimes_RM}=0$ and $\hh{\Hom RKM}$, so the same holds for any element in ${\fa}$.

\begin{noname}\label{koszul amplitude}
Let $(R,\fm,k)$ be a local ring, $M$ a complex of $R$-modules, and let $K$ be the
Koszul complex on a set of $n$ elements in $\fm$. If $\,\hh M$ is degreewise finite,
then
\begin{gather*}
\sup M \leq\sup\left(K\otimes_RM\right)\leq \sup M + n  \,;\tag{1} \\
\inf \left(K\otimes_RM\right) = \inf M\,.\tag{2}
\end{gather*}
\end{noname}

For any element $x$ in $R$, it follows from \eqref{les:koszul} that since $\hh M$ is
degreewise finite, so is $\hh{K(x)\otimes M}$. So it suffices to verify the claim for
$K=K(x)$. Then, invoking \eqref{les:koszul} once again, one deduces that
$\sup(K\otimes_RM)\leq \sup M + 1$ and $\inf \left(K\otimes_RM\right)\geq\inf M$. By
Nakayama's lemma, since $\hh M$ is degreewise finite, $\HH n{K\otimes_RM}=0$ entails
$\HH nM=0$. This yields the remaining inequalities.

\subsection*{The derived category}
As stated in the introduction, our arguments usually take place in $\dcat R$, the
derived category of $R$-modules. Recall that the objects of $\dcat R$ are complexes of
$R$-modules. We use the symbol $\simeq$ to signal an isomorphism in $\dcat R$.  A
morphism of complexes (in the category of complexes of $R$-modules) represents an
isomorphism in $\dcat R$ precisely when it is a quasiisomorphism.

The left derived functor of the tensor product functor of complexes of $R$-modules is
denoted $\Ltor {-}R{-}$, and the right derived functor of the homomorphisms functor is
denoted $\Rhom R{-}{-}$. Given complexes $M$ and $N$, the complexes $\Ltor MRN$ and
$\Rhom RMN$ are  defined uniquely (up to canonical isomorphism) in $\dcat R$. Set
\[
\Tor iRMN = \HH i{\Ltor MRN} \qnd \Ext iRMN = \HH {-i}{\Rhom RMN}\,.
\]
When $M$ and $N$ are modules, these coincide with the classical notions. The reader can
refer to Hartshorne \cite{Ha}, or Gelfand and Manin \cite{GeM}, for a thorough discussion
of the construction and basic properties of the derived category and derived functors.

Given our context, it is crucial that the derived functors are defined without any
boundedness restrictions imposed on the homology of the complexes involved. This is
feasible thanks to the work of Spaltenstein \cite{Sp}, wherein he constructs
appropriate projective and injective resolutions. For the purpose on hand, we do not
need to be concerned with the precise structure of such resolutions. All that is
required are the following notions.

Let $M$ be a complex of $R$-modules such that $\hh M$ is bounded below. 
Then $M$ admits a bounded below flat resolution, that is, a complex $F$ isomorphic to $M$
(in $\dcat R$) with each module $F_\ell$ flat, and $F_\ell=0$ for $\ell\ll0$. The
{\em flat dimension} of $M$ is the number
\[
\fdim_RM = \inf\{\,\sup\{\,\ell\mid F_\ell\ne0\,\} \mid F\,\,\text{a flat
resolution of $M$}\,\}\,.
\]
As is well known, when $M$ is a finitely generated module over a local ring $(R,\fm,k)$,
one has that $\fdim_RM=\sup(\Ltor kRM)$; this is a direct consequence of $M$ having a
minimal free resolution. We require the following extension, obtained by Avramov and Foxby
\cite[(5.5.F)]{AF}, of this result.

\begin{chunk}\label{fdim}
{\em Let $\vf\col (R,\fm,k) \to S$ be a local homomorphism and $M$ a complex of $S$-modules
such that $\hh M$ is finite. Then} 
\[
\fdim_RM = \sup(\Ltor kRM)\,.
\]
\end{chunk}

An analogous notion is that of a {\em projective resolution} of $M$, and also of the {\em
projective dimension} of $M$, which we denote $\fdim_RM$.

\begin{chunk}\label{accounts}{\bf Accounting principles.}
Let $M$ and $N$ be complexes of $R$-modules. The following assertions may be verified
without much ado; alternatively, look into \cite{Fo2}.
\begin{equation*}
 \inf(\Ltor MRN) \geq \inf M + \inf N\,;\tag{1}
\end{equation*}
equality holds if $i=\inf M$ and $j=\inf N$ are finite, and $\HH iM\otimes_R\HH
jN\ne0$.
\begin{equation*}
 \sup(\Ltor MRN) \leq \sup M + \fdim_RN\,;\tag{2}
\end{equation*}
\begin{equation*}
 \sup\Rhom RMN \leq \sup N - \inf M\,;\tag{3}
\end{equation*}
equality holds $i=\inf M$ and $s=\sup N$ are finite and $\Hom R{\HH iM}{\HH sN}\ne 0$.

\noindent{(4)} If $\hh M$ is bounded, then
\[
\sup \Rhom RMN \leq \sup\{\,\sup\Rhom R{\HH {\ell}M}N - \ell\mid \ell\in\BZ\,\}\,.
\]
 \end{chunk}
 
 For ease of reference, we now recall a few facts concerning local cohomology and derived
 completions. Some of these are well known for modules, and perhaps also for bounded
 complexes. In the generality we desire, they are to be found in the paper of Alonso,
 Jerem{\'{\i}}as, and Lipman \cite{Li2}. Having said this, in what follows we refer to Lipman's
 article \cite{Li} whenever possible.

\subsection*{Local cohomology}
Let ${\fa}$ be an ideal in $R$ and $M$ an $R$-module. The {\em ${\fa}$-torsion
submodule of $M$} is the $R$-module
\[
\Gamma_{\fa}(M) = \{\,m\in M\mid {\fa}^n\5 m=0\quad \text{for some integer $n$}\,\} \] The
association $M\mapsto \Gamma_{\fa}(M)$ extends to define an additive, left exact functor
on the category of complexes of $R$-modules; its right derived functor is denoted $\Lch
{\fa}-$. More precisely, given a complex of $R$-modules $M$, let $I$ be an appropriate
injective resolution of $M$ and set $\Lch{\fa}M=\Gamma_{\fa}(I)$. There is a natural morphism
\[
 \Lch {\fa}M \to M \quad \text{in} \quad \dcat R\,.
\]
It is traditional to set
\[
\lch i{\fa}M = \HH {-i}{\Lch {\fa}M}\, \qquad \text{for each}\quad i\in\BZ\,.
\]
This is the {\em $i$th local cohomology $M$ with support in ${\fa}$}. We need the
following properties of local cohomology.

\begin{noname}\label{identity:gamma}
If \,${\fa}\5 \hh M=0$ then the morphism $\Lch {\fa}M\to M$ is an isomorphism.
\end{noname}
This is a special case of \cite[(3.2.1)]{Li}. The argument runs as follows: One
has a spectral sequence that lies in a rectangular strip in the left half-plane, with
\[
\EH 2pq = \lch {-p}{\fa}{\HH qM} \qnd \ED rpq\col \EH rpq \to \EH r{p-r}{q+r-1}\,,
\]
and converges to $\lch {-p-q}{\fa}M$. For each $q$, since ${\fa}\5 \HH qM=0$, one can
construct an injective resolution $I$ of $\HH qM$ with $\Gamma_{\fa}(I)=I$. This yields
$\lch 0{\fa}{\HH qM}=\HH qM$ and $\lch {-p}{\fa}{\HH qM}=0$ for $p\ne0$.  Thus, the
spectral sequence collapses and the edge homomorphisms $\lch{-q}{\fa}M \to \HH qM$ are
isomorphisms.

\begin{noname}\label{commutes:gamma}
For any complex of $R$-modules $N$ there is a canonical morphism
 \[
 \Lch {\fa}{\Rhom RMN} \to \Rhom RM{\Lch {\fa}N}\,;
 \]
it is an isomorphism\, when $\hh M$ is finite and $\pdim_RM<\infty$.
\end{noname}
Indeed, one has an inclusion $\Gamma_\fa \Hom RMN \subseteq \Hom RM{\Gamma_\fa N}$, as is
immediate from the definition of the torsion submodule. Its derived version is the
morphism given above. It is evidently an isomorphism when $M=R$ and hence also when $M$ is
a finitely generated projective module. When both $\hh M$ and $\pdim_RM$ are finite, $M$
can be replaced by a bounded complex of such modules; this explains the assertion about
the isomorphism. 

The displayed morphism is an isomorphism also when $\hh M$ is degreewise finite and
bounded below and $\hh N$ is bounded above, and, in this generality, it
contains Grothendieck's local duality theorem; see \cite[(6.1), (6.3)]{Fo1}. A sheaf theoretic
analogue of the isomorphism above is given by \cite[(5.2.1)]{Li2}. With no restrictions on
$M$ or on $N$, one cannot expect the morphism in \eqref{commutes:gamma} to be an
isomorphism.

We provide two examples to substantiate our claim; both involve a complete local ring $R$,
with maximal ideal $\fa$, and the injective hull $E$ of its residue field. For the first,
we let $M=E$ and $N=E$, so that $\Lch{\fa}N=N$. Thus, the complexes $\Rhom RMN$ and $\Rhom
RM{\Lch{\fa}N}$ are both isomorphic to $R$, and the morphism in \eqref{commutes:gamma} is
the canonical morphism $\Lch{\fa} R \to R$. This last morphism is not an isomorphism
unless $R$ happens to be zero dimensional. This example shows that one
cannot do away with the degreewise finiteness of $\hh M$.

Our second example illustrates the role of the boundedness hypothesis on $\hh N$ when $\hh
M$ is not bounded. Let $M$ and $N$ be complexes with trivial differentials, with
$M_n=R$ and $N_n=E$ for each non-negative integer $n$ and $0$ otherwise. Now, $\Lch{\fa}N=N$
and $M$ is a bounded below complex of projectives, so $\Rhom RMN=\Hom RMN$. This last
complex has trivial differentials, so it is easy to compute its homology:
\[
\HH n{\Hom RMN}=\Hom RMN_n =\prod_{i\ges n}E \quad\text{for each integer $n$.}
\]
In particular, the modules $\HH n{\Hom RMN}$ are not $\fa$-torsion, unless $R$ happens to
be artinian. However the homology modules of $\Lch{\fa}{\Rhom RMN}$ are all $\fa$-torsion;
this is a particular case of the following phenomenon: for any complex $X$, each local
cohomology module $\lch i{\fa}X$ is $\fa$-torsion. This result is contained in the
work of Dwyer and Greenlees; cf.~\cite[(5.3)]{DG}.

\subsection*{Derived completions}
Let ${\fa}$ be an ideal in $R$ and $M$ an $R$-module. The {\em ${\fa}$-adic completion of
  $M$} is the $R$-module
\[
\Lambda_{\fa}(M) = \varprojlim_n (M\otimes_RR/{\fa}^n)\,.
\]
The mapping $M\mapsto\Lambda_{\fa}(M)$ extends to an additive functor on the category of
complexes of $R$-modules. This functor admits a left derived functor that we denote
$\Lho{\fa}-$ following \cite{Li}. This is defined as follows: given a complex of
$R$-modules $M$, let $F$ be an appropriate flat resolution of $M$, and set
$\Lho{\fa}M=\Lambda_{\fa}(F)$.

It is not entirely obvious that this construction yields a well defined functor: given two
flat resolutions $F$ and $G$ of $M$, one has to prove that $\Lambda_{\fa}(F)\simeq
\Lambda_{\fa}(G)$. This entails proving that for any complex of $R$-module $X$, if $\hh
X=0$, then $\hh{\Lambda_{\fa}(X)}=0$ as well. The latter assertion may be verified, for
example, by arguments akin to those used in the last part of the proof of Theorem
\ref{equality:width}. There is another option: one way to define the left derived functor of
$\Lambda_{\fa}$, and indeed, any additive functor on the category of complexes, is via
projective resolutions. Then, since any two such are homotopy equivalent, it is evident
that the functors thus obtained are well defined. However, for most applications it is
crucial that one be able to compute the left derived functor of $\Lambda_{\fa}(-)$ via
flat resolutions.

For each complex of $R$-modules $M$, there is a natural morphism
\[
  M \to \Lho {\fa}M \qquad\text{in}\quad \dcat R\,.
 \]
For each integer $i$, the {\em $i$th derived completion of $M$ with respect to ${\fa}$}
is the $R$-module
\[
\lho i{\fa}M = \HH i{\Lho {\fa}M}\,.
\]
We require the following facts concerning derived completions.

\begin{noname}\label{identity:lambda}
If \,${\fa}\5\hh M=0$ then the morphism $M\to \Lho {\fa}M$ is an isomorphism.
\end{noname}

Here is one justification of this assertion: By \cite[(0.3)]{Li2}, see also
\cite[(4.1)]{Li}, there is a natural isomorphism $\Rhom R{\Lch {\fa}R}M\to \Lho
{\fa}M$. The projective dimension of $\Lch {\fa}R$ is finite;
cf.~\cite[(6.5)]{Fo1}. Thus, a spectral sequence analogous to the one in
\eqref{identity:gamma} allows one to reduce the problem to the case where $M$ is
an $R$-module with ${\fa}\5 M=0$. Then the canonical map $M\to \Lambda_{\fa}M$
is an isomorphism so the same holds for the morphism $M\to\Lho {\fa}M$, by
\cite[(4.1)]{GM}.

Theorem \ref{identity:lambda} is valid under the far weaker assumption that for each
integer $i$ the $R$-module $\HH iM$ is $\fa$-adically complete. Similarly, Theorem
\ref{identity:gamma} holds, more generally, whenever each $\HH iM$ is $\fa$-torsion.  A
good way to understand these phenomenon is via the theory of Bousfield colocalizations and
localizations; confer, for example, the article of Dwyer and Greenlees \cite{DG}.

\begin{noname}\label{NAK:lambda}
For any integer $i$, if ${\fa}\5\lho i{\fa}M=\lho i{\fa}M$, then $\lho
i{\fa}M=0$.
\end{noname}

Indeed, this follows from \cite[(1.4)]{Si} as in \cite{Fr}, that is, by replacing
$M$ by an appropriate resolution $P$ and noting that then $\lho i{\fa}M$ equals $\HH
i{\Lambda_{\fa}(P)}$.


\begin{noname}\label{commutes:lambda}
For any complex of $R$-modules $N$, there is a canonical morphism
\[
 \Ltor{\Lho {\fa}M}RN\to\Lho {\fa}{\Ltor MRN}\,;
 \]
it is an isomorphism, if $\hh N$ is finite and $\pdim_RN<\infty$.
\end{noname}

As noted in \eqref{identity:lambda} above, the complexes $\Rhom R{\Lch {\fa}R}M$ and
$\Lho {\fa}M$ are canonically isomorphic. The desired morphism results from the natural
morphism
\[
\Ltor{\Rhom R{\Lch {\fa}R}M}RN \to \Rhom R{\Lch {\fa}R}{\Ltor MRN}\,.
\]
As to the bit about the isomorphism: It is immediate when $N=R$, and hence also when
$N$ is a finitely generated projective module. This settles it, since $N$ is isomorphic
in $\dcat R$ to a bounded complex consisting of such modules.

\section{Depth}\label{depth}

The following is one of our main results. It contains Theorem I discussed in the
introduction, and extends results of Iyengar \cite[\S6]{Iy} who operates under the
additional hypothesis that $\hh M$ is bounded above.

\begin{theorem}\label{equality:depth}
Let $R$ be a noetherian ring and $M$ a complex of $R$-modules. Let ${\fa}$ be an ideal
in $R$ and let $K$ be the Koszul complex on a sequence of $n$ generators for ${\fa}$. In
this case, one has that
\begin{align*}
\sup(K\otimes_RM)-n &= \sup\Hom RKM \\
                    &= \sup\Rhom R{R/{\fa}}M \\
                    &= \sup\Lch {\fa}M\,.
\end{align*}
\end{theorem}

The proof of this theorem uses the proposition below.

In fact, this latter result can be extended so that it holds, more generally, for any
$X$ such that $\Lch {\fa}X\simeq X$, but justifying this last claim
requires the use of sophisticated tools.  The same comment applies also to \eqref{weak
sensitivity:width}. What is more, these results, and others in their vein, have
implications that go beyond the present application; we plan to turn to these matters
in future work.

For now, we state and prove only the weaker version below for it suffices for the
present purpose and its proof is elementary.

\begin{proposition}\label{weak sensitivity:depth}
Let $M$ and $X$ be complexes of $R$-modules. If $\,\hh X$ is bounded and ${\fa}^d\5 \hh
X=0$ for some positive integer $d$, then
\[
\sup \Rhom RXM \leq \sup \Rhom R{R/{\fa}}M - \inf X\,.
 \]
\end{proposition}

\begin{proof}
Since $\hh X$ is bounded, (\ref{accounts}.4) yields the estimate
\[
\sup \Rhom RXM\leq \sup\{\,\sup\Rhom R{\HH {\ell}X}M- \ell\mid \ell\in\BZ\,\}\,.
\]
This reduces the problem to the case where $X$ is concentrated in degree $0$. By
hypothesis ${\fa}^d\5 X=0$, so $X$ is in fact an $R/{\fa}^d$-module. For any
$R/{\fa}^d$-module $T$, one has the isomorphism
\[
 \Rhom RTM \simeq \Rhom {R/{\fa}^d}T{\Rhom R{R/{\fa}^d}M}\,.
\]
By (\ref{accounts}.3), this implies that
\begin{equation*}
\sup\Rhom RTM \leq \sup \Rhom R{R/{\fa}^d}M \,.\tag{$\dagger$}
 \end{equation*}
Therefore, it suffices to prove the result in the case where $X=R/{\fa}^d$. To do this,
we resort to an induction on the integer $d$, the base case $d=1$ being tautological.
Suppose that the desired estimate holds for $X=R/{\fa}^d$ for some integer $d\geq 1$ .
Applying $\Hom R{-}M$ to the short exact sequence
\[
0\to {{\fa}^d}/{{\fa}^{d+1}}\to R/{{\fa}^{d+1}}\to R/{{\fa}^d} \to 0\,,
 \]
yields the long exact sequence
\[
\cdots\to\Ext mR{R/{\fa}^d}M\to\Ext mR{R/{\fa}^{d+1}}M\to \Ext
mR{{\fa}^d/{\fa}^{d+1}}M\to \cdots\,.
\]
Since ${\fa}^d/{\fa}^{d+1}$ is an $R/{\fa}^d$-module, ($\dagger$) yields that
$\sup\Rhom R{{\fa}^d/{\fa}^{d+1}}M \leq \sup\Rhom R{R/{\fa}^d}M$. This estimate, along
with the induction hypothesis and the long exact sequence above allows us to complete
the induction step, and hence the proof of the proposition.
\end{proof}

\begin{proof}[Proof of Theorem \ref{equality:depth}]
Since $K$ is a bounded complex of finitely generated free $R$-modules, the complexes
$\Ltor KRM$ and $K\otimes_RM$ are isomorphic. For this reason, in the ensuing
discussion, we identity these two complexes; ditto for the complexes $\Rhom RKM$ and
$\Hom RKM$.

{\em Proof of\/} $\sup(K\otimes_RM)-n=\sup\Hom RKM$.  Koszul complexes are self dual:
There is an isomorphism $K\cong \shift n{\Hom RKR}$ of complexes of $R$-modules and
this induces the isomorphism $K\otimes_RM\cong \shift n{\Hom RKM}$. In particular,
\[
\sup\left(K\otimes_RM\right) = \sup \Hom RKM + n\,.
\]

{\em Proof of\/} $\sup\Hom RKM = \sup\Rhom R{R/{\fa}}M$. Since $\hh K$ is bounded with
${\fa}\5\hh K=0$, the proposition above yields the inequality
\[
\sup\Rhom RKM \leq \sup \Rhom R{R/{\fa}}M\,.
\]
As to the opposite inequality, note that
\begin{align*}
\sup\Rhom R{R/{\fa}}{\Rhom RKM} &=\sup \Rhom R{\Ltor{(R/{\fa})}RK}M \\
                            &=\sup\Rhom R{\coprod_{i\ges 0}\shift
i{(R/{\fa})^{\binom ni}}}M\\
                            &=\sup \Rhom R{R/{\fa}}M\,.
\end{align*}
The desired result follows from the calculation above and the fact that, by
(\ref{accounts}.3), $\sup\Rhom R{R/{\fa}}{\Rhom RKM}\leq \sup\Rhom RKM$.

{\em Proof of\/} $\sup \Hom RKM = \sup \Lch {\fa}M$.  In the chain below, the
isomorphism on the left is by \eqref{identity:gamma}, since ${\fa}\5 \hh{\Hom RKM}=0$,
while the one on the right is given by \eqref{commutes:gamma}.
\[
\Hom RKM \simeq \Lch {\fa}{\Hom RKM}\simeq\Rhom RK{\Lch {\fa}M}\,.
\]
This yields the equality below, while the inequality follows from (\ref{accounts}.3).
\[
\sup\Hom RKM = \sup \Rhom RK{\Lch {\fa}M} \leq \sup\Lch {\fa}M\,.
 \]
On the other hand, for each integer $n\geq 0$, one has that
\[
\sup\Rhom R{R/{\fa}^n}M \leq \sup \Rhom R{R/{\fa}}M = \sup \Hom RKM\,,
\]
where the equality has been established already, and the inequality is provided by
\eqref{weak sensitivity:depth}.  Since $\lch i{\fa}M = \varinjlim{}_n\Ext iR{R/{\fa}^n}M$
for each integer $i$, the inequalities above imply that $\sup\Lch {\fa}M \leq \sup\Hom
RKM$.

This completes the proof of the desired equality and of the theorem.
 \end{proof}

\subsection*{A compendium}
As has been explained in the introduction, Theorem \ref{equality:depth} may be
interpreted as stating that the various ways of introducing depth for complexes: via
Koszul homology and cohomology, via the {\rm Ext} functor, and via local cohomology,
all lead to the same invariant. Thus, there is no ambiguity in speaking of {\em the}
depth of a complex. Let us record this fact.

\begin{definition}\label{definition:depth}
Let $R$ be a noetherian ring, ${\fa}$ an ideal in $R$ and $K$ the Koszul complex on a
sequence of $n$ generators for ${\fa}$. For a complex of $R$-modules $M$, the {\em
${\fa}$-depth of $M$ over $R$} is defined by one of the following equivalent formulas:
\[
\Rdepth R{\fa}M = \begin{cases}
   n - \sup(K\otimes_RM)\,; \cr
  \inf\{\,\ell\in\BZ\mid \HH{-\ell}{\Hom RKM}\ne0\,\}\,; \cr
  \inf\{\,\ell\in\BZ\mid \Ext \ell R{R/{\fa}}M\ne 0\,\}\,; \cr
    \inf\{\,\ell\in\BZ\mid \lch \ell {\fa}M\ne 0\,\}\,.
\end{cases}
\]
We write $\depth {\fa}M$ for the ${\fa}$-depth of $M$ when the ring $R$ is clear from
the context. By the by, the equality above implies that in computing depth via Koszul
(co)homology, one may choose any finite generating sequence for ${\fa}$.

For a complex $M$ over a local ring $(R,\fm,k)$, the {\em depth of $M$} is the number
\[
\mdepth_R M =\Rdepth R{\fm}M\,.
\]
This is abbreviated to $\mdepth M$, if omitting $R$ does not lead to much confusion.
\end{definition}

In the remainder of this section we state, and prove, the fundamental properties that
depth enjoys without imposing unnecessary boundedness conditions. In this process, it
becomes clear that in dealing with depth no one definition can be singled out as being
best suited to every purpose.

The first result subsumes Theorem II from the introduction. When $\hh M$ is
bounded above, it is precisely \cite[(2.1)]{Iy}; in turn that extends results in
\cite{Fo2}, \cite{Iv2}.

\begin{theorem}\label{AB:depth}
Let $(R,\fm,k)$ be a local ring and $M$ a complex of $R$-modules. Let $P$ be bounded
complex of $R$-modules with $\hh P\ne0$ and such that $\fdim_RP<\infty$.

If either $\hh M$ is bounded above or $\hh P$ is degreewise finite, then
\[
\mdepth_R(\Ltor MRP) = \mdepth_RM - \sup(\Ltor kRP) \,.
\]
 \end{theorem}
\begin{proof}
  As noted above, the case where $\hh M$ is bounded above is settled by \cite[(2.1)]{Iy};
  confer also \cite[(1.8)]{Fo1}.  That proof works with the Koszul homology
  characterization of depth. We give an argument that handles both parts simultaneously.

The complex of $R$-modules $\Rhom RkM$ is isomorphic to a graded $k$-vector space; for
example, see \cite{Fo2}. This engenders the isomorphism in the diagram below, whereas
the morphism $\theta$ is the canonical one.
\[
\Ltor{\Rhom RkM}k{(\Ltor kRP)}\simeq \Ltor{\Rhom RkM}RP
               \arto{\theta} \Rhom Rk{\Ltor MRP}\,.
\]
Under either hypothesis, $\theta$ is an isomorphism: When $\hh M$ is bounded above,
this is \cite[(4.4.F)]{AF}; a similar argument also goes through when $\hh P$ is
degreewise finite.  If $\hh M$ is bounded above, then $\sup\Rhom RkM<\infty$, by
(\ref{accounts}.3). When $\hh P$ is degreewise finite, $\hh{\Ltor kRP}$ is non-zero, by
(\ref{accounts}.1), so (\ref{accounts}.2) yields that $\sup(\Ltor kRP)$ is finite. Thus,
the isomorphisms above implies an equality
\[
 \sup\Rhom Rk{\Ltor MRP} = \sup\Rhom RkM + \sup(\Ltor kRP)\,.
 \]
This is the result we seek.
\end{proof}

At this point, it is expedient to record the following remarks which are handy for many
of the subsequent arguments.

\begin{noname}\label{infinity:depth}
Let $R$ be a local ring, ${\fa}$ a proper ideal in $R$, and let $M$ be a complex
of $R$-modules with $\hh M$ degreewise finite. Then
\begin{alignat*}{2}
&\Rdepth R{\fa}M = -\infty& &\iff \sup M = \infty\,;\\
&\Rdepth R{\fa}M = \infty &  & \iff \sup M = -\infty\,.
 \end{alignat*}
\end{noname}
These  are immediate from (\ref{koszul amplitude}.1), once we compute depth via Koszul
homology.

It is crucial that $R$ be local, as the following example illustrates.

\begin{example}\label{example:depth} Let $R$ be a noetherian ring containing a set
of non-trivial ideals $\{\,{\fa}_n\,\}_{n\ges 0}$ with the property that
${\fa}_i+{\fa}_j=R$ for $i\ne j$, and set
\[
M=\coprod_{n\ges 0} \shift n{(R/{\fa}_n)}\,.
\]
Then, $\hh M$ is degreewise finite with $\sup M=\infty$, while for any non-negative
integer $n$, one finds that $\Rdepth R{{\fa}_n}M = -n$ (compute via Koszul complexes).
\end{example}

The following lower bound for depth is well known; see \cite{Fo2} or \cite[(2.3)]{Iy}.

\begin{noname}\label{lower bound:depth}
Let $(R,\fm,k)$ be a local ring and $M$ a complex such that $\sup M =s$ is finite.
Then $\mdepth M\geq -s$ and equality holds if and only if $\fm\in\ass\HH sM$.
\end{noname}

Here is what we have to say about upper bounds for depth.

\begin{chunk}\label{upper bound:depth}
Let ${\fa}$ be an ideal in $R$ and set $a=\max\{\,i\mid \lch j{\fa}R\ne0\,\}$. Let $M$
a complex of $R$-modules. From the local cohomology spectral sequence encountered in
\eqref{identity:gamma} it follows that if there are integers $s$ and $d$ such that
\begin{enumerate}[\quad\rm(a)]
\item $\lch d{\fa}{\HH sM}\ne 0$, whilst
\item $\lch {d-1-j}{\fa}{\HH{s-i}M}=0$ for $1 \leq j \leq d-1 $, and
\item $\lch {d+1+j}{\fa}{\HH{s+j}M}=0$ for $1 \leq j \leq a-d-1$
\end{enumerate}
then $\EH 2{-d}s$ survives to $\EH {\infty}{-d}s$ so $\lch {d-s}{\fa}M\ne 0$. Thus,
$\Rdepth R{\fa}M \leq d-s$.

Conditions (a)--(c) may seem contrived, but they hold in the following case.

\begin{Lemma}
Let $(R,\fm,k)$ be a local ring and $M$ be a complex of $R$-modules. If there is an
integer $s$ such that
 \begin{enumerate}[{\quad\rm(a)}]
\item $\fm\in\ass\HH sM$\ and
\item $\HH jM$ is $\fm$-torsion for each integer $s+1\leq j\leq s+\dim
R - 1$\,,
\end{enumerate}
then $\mdepth_RM \leq -s$. In particular, if $\sup M=\infty$ and $\HH jM$ is
$\fm$-torsion for all $j\gg0$, then $\mdepth_RM=-\infty$. \qed
\end{Lemma}

\end{chunk}

Next we improve on \cite[(5.3)]{Iy}. The proof of \emph{loc.~cit.} uses the
Auslander-Buchsbaum equality \cite[(2.1)]{Iy}, which is why its validity was restricted to
complexes with bounded above homology. Now, thanks to \eqref{AB:depth}, the same
argument establishes this result without any boundedness hypothesis on the
complex $M$.

\begin{proposition}\label{last homology}
Let ${\fa}$ be an ideal in a noetherian ring $R$, and let $K$ be the Koszul complex on
a finite sequence $\bsx$ of generators for ${\fa}$. Let $M$ be a complex of $R$-modules such
that its ${\fa}$-depth is finite. Set $d=\Rdepth R{\fa}M$ and $s=\sup(K\otimes_RM)$.
\begin{enumerate}[{\quad\rm(1)}]
\item The $R$-module $\HH s{K\otimes_RM}$ is independent of the choice
of $\bsx$.
\item $\{\,\fp\in\V {\fa}\mid\mdepth_{R_\fp}M_{\fp}=d\,\}=\ass_R\HH
s{K\otimes_RM}$.\qed
\end{enumerate}
\end{proposition}

The proposition below describes the local nature of depth. It builds on the proof of
\cite[(5.4)]{Iy} which is good enough to handle the case when $\Rdepth R{\fa}M >
-\infty$. Thus the only situation that remains to be tackled is when $\Rdepth
R{\fa}M=-\infty$. This is more involved than one might suspect.

\begin{proposition}\label{locality:depth}
Let ${\fa}$ be an ideal in a noetherian ring $R$ and let $M$ be a complex of
$R$-modules. Then
\[
 \Rdepth RIM = \inf\{\mdepth_{R_\fp}M_{\fp}\mid \fp\in\V I\}\,.
 \]
\end{proposition}

\begin{Remark}
If ${\fa}=0$, then $\Rdepth R{\fa}M=-\sup M$, and the proposition reads
\[
-\sup M = \inf\{\,\mdepth_{R_\fp}M_{\fp}\mid \fp\in\spec R\,\}\,.
\]
This formula in conjunction with the preceding proposition explains \eqref{lower
bound:depth}.
\end{Remark}

\begin{caveat}
Evidently, when $\Rdepth R{\fa}M$ is finite, the infimum is achieved at some prime
$\fp\in \V {\fa}$.  This {\em need not be the case if $\Rdepth R{\fa}M=-\infty$}, as is
illustrated by the following example.

Let $R$ and $M$ be as in \eqref{example:depth}. Then $\depth{(0)}M=\sup M = -\infty$,
whilst $\hh{M_\fp}$ is bounded above for any prime ideal $\fp\in\spec R$, so
$\mdepth_{R_\fp}M_\fp>-\infty$. Note that if $\fp_i$ is a minimal prime in
$\V{{\fa}_i}$, then $\mdepth_{R_{\fp_i}}M_{\fp_i}=-i$, so the infimum here is
$-\infty$, as predicted by the proposition.

It is not too hard to cook up another such example wherein the ring $R$ is {\em local}.
 \end{caveat}

\begin{proof}[Proof of Proposition \ref{locality:depth}]
As mentioned before, \cite[(5.4)]{Iy} resolves the case when $\Rdepth R{\fa}M >
-\infty$.  The argument is short, so it bears repeating: Utilizing the
characterization of depth via either local cohomology or Koszul homology it is easy to
establish the following inequalities; cf.~\cite[(5.2)]{Iy}.
\[
\Rdepth R{\fa}M \leq \Rdepth{R_\fp}{{\fa}_\fp}{M_\fp}\leq \mdepth_{R_\fp}M_\fp\,.
\]
Now we may assume that $\Rdepth R{\fa}M$ is finite, in which case Proposition \ref{last
homology} provides us with a prime $\fp$ at which the inequalities above become
equalities.

For the remainder of the proof $\Rdepth R{\fa}M=-\infty$; in particular, $\sup
M=\infty$. We have to establish that $\mdepth_{R_\fp}M_\fp$ has no lower bound as $\fp$
varies in $\V {\fa}$. To this end, it is convenient to engineer ourselves into a
situation where $M$ is supported in $\V {\fa}$, by the following device: Let $K$ be the
Koszul complex on a sequence of $n$ generators for ${\fa}$. Then
\begin{align*}
 \depth {\fa}{K\otimes_RM}&=n-\sup(K\otimes_R(K\otimes_RM))\\
                      &=n-(2n-\depth {\fa}M)\\
                      &=\depth {\fa}M-n\,,
\end{align*}
where the first, respectively, the second, equality is due to the fact that $K$,
respectively, $K\otimes_RK$, detects depth with respect to ${\fa}$. Thus, $\depth
{\fa}{K\otimes_RM}=-\infty$ as well. Moreover, for any prime $\fp\in \V {\fa}$, the
complex $K_\fp$ is minimal, so that $\sup\left(\Ltor{k(\fp)}{R_\fp}{K_\fp}\right)=n$,
where $k(\fp)$ denotes the residue field of $R_\fp$. In particular, Theorem
\ref{AB:depth} yields
\[
\mdepth_{R_\fp}(K_\fp\otimes_{R_\fp}M_\fp)= \mdepth_{R_\fp}M_\fp - n\,.
 \]
Therefore, by passing to $K\otimes_RM$ we may assume in addition that ${\fa}\5\hh M=0$.

Define $V_n=\{\,\fq\in\V {\fa}\mid n\leq \sup M_\fq < \infty\,\}$, for each integer
$n$. These are subsets of $\spec R$ with $\cdots\supseteq V_n\supseteq V_{n+1}\supseteq
\cdots$. There are two possibilities.

Suppose that $V_n\ne\varnothing$ for each integer $n$. Fix an integer $d\geq 0$ and
choose an element $\fq\in V_d$, so that $\sup M_\fq=s \geq d$ and finite. If $\fp$ is a
prime associated to $\HH s{M_\fq}$, then ${\fa}\subseteq \fp$ because ${\fa}\5\hh M=0$.
Moreover, we deduce from \eqref{lower bound:depth} that $\mdepth_{R_\fp}M_\fp=-s\leq -d$,
which yields the desired conclusion, since $d$ was arbitrary.

Suppose that $V_n=\varnothing$ for some integer $n$. Then, since ${\fa}\5\hh M=0$ and
$\sup M=\infty$, this implies that the set $U=\{\,\fq\in\V {\fa}\mid \sup M_\fq
=\infty\,\}$ is non-empty. Pick an prime ideal $\fp$ which is minimal in $U$, that is
to say, $\fq\not\subset\fp$ for any element $\fq\in U$. By choice of $\fp$, one has
$\sup M_\fp = \infty$, whilst for $i\geq n$ the $R_\fp$-module $\HH i{M_\fp}$ is
supported only at the maximal ideal $\fp R_\fp$. It remains to invoke Lemma \ref{upper
bound:depth} to conclude that $\mdepth_{R_\fp}M_{\fp}=-\infty$.

This completes the proof of the proposition.
\end{proof}

The preceding proposition allows us to extend \cite[(5.5)]{Iy} to complexes which are
not necessarily bounded above.

\begin{proposition}\label{depth:compare}
Let ${\fa}$ and ${\fb}$ be ideals in a noetherian ring $R$, and let $M$ be a complex of
$R$-modules.
\begin{enumerate}[{\quad\rm(1)}]
\item $\depth {{\fa}{\fb}}M =\depth{{\fa}\cap {\fb}}M=\min\{\,\depth
{\fa}M,\depth {\fb}M\,\}$.
\item $\depth {\fb}M = \depth {\fa}M$ if $\,\rad({\fb})=\rad({\fa})$.
\end{enumerate}

If $R$ is local and $\hh M$ is degreewise finite, then
\begin{enumerate}[{\quad\rm(1)}]
\item[\rm(3)] $\mdepth M \leq \depth {\fa}M + \dim R/{\fa}\,.$
\end{enumerate}
\end{proposition}
\begin{proof}
(1) This first equality is a corollary of \eqref{locality:depth} since
$\V{{\fa}{\fb}}=\V{{\fa}\cap {\fb}}$. As to the second, consider the equalities
\begin{align*}
\depth{{\fa}\cap {\fb}}M&=
\inf_{\fp\in\V{{\fa}\cap{\fb}}}\mdepth_{R_\fp}M_\fp\\
 &= \min\{\,
\inf_{\fp\in\V{{\fa}}}\mdepth_{R_\fp}M_\fp\,,\,
\inf_{\fp\in\V{{\fb}}}\mdepth_{R_\fp}M_\fp\,\}\\
&=\min\{\,\depth {\fa}M,\depth {\fb}M\,\}\,,
\end{align*}
where the first and the third are by \eqref{locality:depth}, and the middle one reflects
the identity: $\V{{\fa}\cap {\fb}}=\V {\fa} \cup \V {\fb}$.

(2) This too follows from \eqref{locality:depth} since $\V {\fa}=\V {\fb}$.

(3) Thanks to \eqref{infinity:depth}, it suffices to consider the case when $\sup M$ is
finite. At this point, we may refer to \cite[(5.5.4)]{Iy}, but for completeness we give (a
slight variant of) the argument: Since $R$ is local, one can find elements $x_1,\dots,x_d$
in $R$ whose image under the canonical surjection $R\to R/{\fa}$ forms a system of
parameters for $R/{\fa}$. In particular, $\rad(\bsx,{\fa})$ is the maximal ideal $\fm$ of
$R$, and $d=\dim R/{\fa}$. Thus, $\Lch {\fm}M \simeq\Lch {(\bsx)}{\Lch {\fa}M}$; for
example, see \cite
{Li}. Since $\sup M<\infty$, one has $\sup\Lch {\fa}M<\infty$, so $\sup
\Lch{(\bsx)}{\Lch {\fa}M}\geq \sup \Lch {\fa}M - d$, which is the inequality we seek.
\end{proof}

\section{Amplitude inequality}\label{amplitude}

This section is dedicated to the proof of the following extension of Theorem III from
the introduction.

\begin{theorem}\label{amplitude inequality}
Let $R\to S$ be a local homomorphism and let $F$ be a complex of $S$-modules with $\hh
F$ non-trivial and finite. If the flat dimension of $F$ over $R$ is finite, then for
any complex of $R$-modules $M$ with $\hh M$ degreewise finite, one has
\begin{gather*}
\sup M+\inf F  \leq  \sup\left(\Ltor MRF\right) \,; \\
 \inf M+\inf F = \inf\left(\Ltor MRF\right)\,.
\end{gather*}
In particular, $\amp M \leq \amp\left(\Ltor MRF\right)$.
\end{theorem}

This theorem has the following surprising (to us) corollary. It contain Theorem IV
stated in the introduction; as has been explained there, this latter result is a
significant generalization of ~\cite[Theorem R]{Ap}.

\begin{theorem}\label{fdim:theorem}
Let $Q\to R\to S$ be local homomorphisms and let $F$ be a complex of $S$-modules such that
$\hh F$ is non-trivial and finite. Assume that $\fdim_RF$ is finite. One has
\[
\fdim_QR + \inf F \leq \fdim_QF \leq \fdim_QR + \fdim_RF\,.
\]
In particular, $\fdim_QR$ and $\fdim_QF$ are finite simultaneously.
 \end{theorem}
\begin{proof}
Let $k$ be the residue field of $Q$, and let $M=\Ltor kQR$. Observe that the complexes
$\Ltor MRF$ and $\Ltor kQF$ are isomorphic. In particular, by \eqref{fdim}, one
has $\fdim_QR=\sup M$ and $\fdim_QF=\sup(\Ltor MRF)$, so the inequality on the left is
a consequence of the preceding theorem. The one on the right is given by
(\ref{accounts}.2).
\end{proof}

  From the preceding result one can deduce (a complex extension of) Theorem V from the
introduction. The equivalence of the first two conditions is contained in a result of Kunz
\cite[(2.1)]{Ku}; that of the first and the third was discovered by Rodicio
\cite[2]{Ro} in the special case where $F=R$. The reader may consult the
survey article of C.~Miller \cite{Mi} in these proceedings for other developments that are
inspired by the work of Kunz and Rodicio.

\begin{theorem}
Let $R$ be a local ring of characteristic $p$, and let $\vf\col R\to R$ be the
Frobenius endomorphism of $R$. The following conditions are equivalent.
\begin{enumerate}[{\quad\rm(a)}]
\item $R$ is regular; 
\item $\vf^n$ is flat for each integer $n\geq 0$;
\item there exists a positive integer $n$ and a complex of $R$-modules $F$ with
$\hh F$ finitely generated such that both $\fdim_RF$ and $\fdim_R\frob nF$ are finite.
\end{enumerate}
\end{theorem}
\begin{proof}
  Let $\fm$ be the maximal ideal of $R$. For each natural number $n$, one has
  $\mdepth_R \frob nR=\Rdepth R{\vf^n(\fm)}R=\mdepth R$. Indeed, the first equality is
  immediate from the Koszul complex characterization of depth whilst the second is given
  by \eqref{depth:compare}, since $\rad(\vf^n(\fm))=\rad(\fm)$. These equalities, coupled
  with equality $\sup(\Ltor kR{\frob nR})=\fdim_R\frob nR$, provided by \eqref{fdim}, and
  Theorem \ref{AB:depth} yield
\begin{equation*}
\fdim(\vf^n) < \infty \implies \fdim(\vf^n)=0\,.\tag{$\dagger$}
\end{equation*}
Now for the proof of the desired equivalences.

(a) $\implies$ (b): Since $R$ is regular, $\fdim\vf^n$ is finite, and hence
$\fdim\vf^n=0$, by ($\dagger$).

(b) $\implies$ (c): Pick a positive integer $n$ and set $F=R$.

(c) $\implies$ (a): Applied to the diagram $R\arto{\vf}R\arto{\vf^{n-1}}R$, Theorem
\ref{fdim:theorem} entails $\fdim(\vf)<\infty$. From this inequality, ($\dagger$) allows us to draw
the stronger conclusion that $\fdim\vf=0$. It remains to invoke \cite[(2.1)]{Ku}.
 \end{proof}

\begin{proof}[Proof of Theorem \ref{amplitude inequality}]
If $\hh M=0$, then $\hh{\Ltor MRF}=0$, and the desired (in)equalities are immediate.
For the rest of the proof it is assumed that $\hh M\ne 0$.

Let $\fm$ be the maximal ideal of $R$ and set $k=R/\fm$. Let $\vf$ denote the
homomorphism $R\to S$. The first step is to reduce to the case where $\vf$ is
surjective; then $\hh F$ would be finite over $R$ itself.

Let $\wh S$ denote the completion of $S$ at its maximal ideal, and set $\wh F =
F\otimes_S{\wh S}$. Since the $S$-module $\wh S$ is faithfully flat, $\hh{\wh F}\cong
\hh F\otimes_S\wh S$; hence $\hh{\wh F}$ is degreewise finite over $\wh S$. For the
same reason, for each complex of $R$-modules $X$, one has
\[
\sup\left(\Ltor XR{\wh F}\right)=\sup\left(\Ltor XRF\right) \qnd \inf\left(\Ltor XR{\wh
F}\right)=\inf\left(\Ltor XRF\right)\,.
\]
Thus, it suffices to establish the (in)equalities we seek with $\wh F$ in place of $F$.
Moreover, the special case $X=k$ of the equality above concerning the suprema, in
conjunction with \eqref{fdim}, yields: $\fdim_R\wh F=\fdim_RF$. Hence the flat
dimension of $\wh F$ over $R$ is finite. So, passing to $\wh S$ and $\wh F$, one can
assume that $S$ is complete.

By \cite[(1.1)]{AFH}, the homomorphism $\vf$ has a factorization
$R\arto{\dot\vf}R'\arto{\vf'}S$ such that the $R$-module $R'$ is flat, the ring $R'/\fm
R'$ is regular, and the map $\vf'$ is surjective.  Set $M'=M\otimes_RR'$.  Then,
$\Ltor{M'}{R'}F\simeq \Ltor MRF$, and, since $\dot\vf$ is faithfully flat, $\sup
M'=\sup M$ and $\inf M'=\inf M$. Furthermore, it follows from arguments analogous to
\cite[(3.2)]{AFH} that the flat dimension of $\fdim_{R'}F$ is finite. Thus, replacing
$R$ and $M$ by $R'$ and $M'$ respectively, one can assume that $\vf$ is a surjective
homomorphism.

 From this point onwards the ring $S$ plays no role in the picture.

The homomorphism $\vf$ is surjective, so $\hh F$ is finite over $R$. Moreover,
$\fdim_RF$ is finite. Thus, $F$ is isomorphic to a complex $F'$ of finitely generated,
free modules with $F'_i=0$ for $i\geq \fdim_RF+1$. Replacing $F$ by $F'$, one can
assume henceforth that $F$ is of this form. In particular, $\Ltor MRF\simeq
M\otimes_RF$.

{\em Proof that $\sup M+\inf F \leq \sup\left(M\otimes_RF\right)$}. If $\sup M=\infty
$, then $\mdepth M = -\infty$ by \eqref{infinity:depth}, so
$\mdepth(M\otimes_RF)=-\infty$, by the Auslander-Buchsbaum formula \eqref{AB:depth};
another application of \eqref{infinity:depth} yields $\sup(M\otimes_RF) =\infty$.  Now we
may assume that $\sup M$ is finite. The next step is to reduce to case where $\inf M$
is also finite, so that $\hh M$ is bounded.

Set $s=\sup(M\otimes_RF)$. We claim that $-\infty < s <\infty $.

Indeed, since $\sup M<\infty$ and $F$ is finite free, $s<\infty$; see (\ref{accounts}.2).
By \eqref{infinity:depth}, the inequality $\sup M > -\infty$ implies $\mdepth M<\infty$,
and hence, by \eqref{AB:depth}, that $\mdepth{(M\otimes_RF)}<\infty$. Another appeal to
\eqref{infinity:depth} yields that $s>-\infty$.

Set $M'=\rtrunc uM$, where $u=s-\pdim_RF - 2$. The following equalities hold.
\begin{enumerate}[{\quad\rm(1)}]
\item $\sup M' = \sup M$;
\item $\sup(M'\otimes_RF) =\sup(M\otimes_RF)$.
\end{enumerate}
To see this, note that the canonical inclusion $M'\to M$ induces the morphism of
complexes $M'\otimes_RF\to M\otimes_RF$. This map is the identity in degrees $s-1$ and
higher so that $\HH n{M'\otimes_RF}=\HH n{M\otimes_RF}$ for $n\geq s$; this contains
$(2)$. In particular, $\hh{M'}\ne 0$ and hence $u\leq \sup M$; thus $\sup{M'}=\sup M$.

At this point, we substitute $M'$ for $M$ and assume that $\hh M$ is bounded. In view
of \cite{Pr}, the result of Iversen \cite[(3.2)]{Iv1} gives
$\amp(M)\leq\amp(M\otimes_RF)$.  Unravelling this inequality, keeping in mind that
$\inf M + \inf F =\inf(M\otimes_RF)$, by (\ref{accounts}.1), yields the desired
inequality.

This completes the justification of the inequality concerning the suprema.

{\em Proof that $\inf M+\inf F = \inf\left(M\otimes_RF\right)$}. When $\inf M$ is
finite, (\ref{accounts}.1) provides the desired equality. Suppose that $\inf M=-\infty$.

Let $K$ be the Koszul complex on a finite generating set for $\fm$. For any complex of
$R$-modules $X$ with $\hh X$ degreewise finite, $\inf\left(K\otimes_RX\right)=\inf X$;
this follows from \eqref{les:koszul}.  In particular,
$\inf\left(K\otimes_RM\right)=-\infty$; moreover,
\[
\inf\left(M\otimes_RF\right)=-\infty \iff
\inf\left(\Ltor{(K\otimes_RM)}RF\right)=-\infty\,.
\]
So, substituting $K\otimes_RM$ for $M$, one can assume that the $R$-module $\HH iM$ has
finite length; see \eqref{annihilation:koszul}. Let $E$ be the injective hull of $k$ -
the residue field of $R$ - and set $M^\vee = \Hom RME$. Then $\HH i{M^\vee}$ has finite
length, and in particular, finitely generated, and $\sup{M^\vee}=-\inf M = \infty$.
This explains the equality in the calculation below; the inequality is by the already
established part of the theorem, since $\Hom RFR$ is finite free.
\[
\sup\left({M^\vee}\otimes_R{\Hom RFR}\right) \geq \sup{M^\vee} + \inf{\Hom
RFR}=\infty\,.
\]
The complexes $M^\vee\otimes_R\Hom RFR$ and $\Hom R{F\otimes_RM}E$ are isomorphic,
since $F$ is finite free. So $\sup \Hom R{F\otimes_RM}E=\infty$; this entails
$\inf\left(M\otimes_RF\right)=-\infty$.

This completes the proof of the theorem.
 \end{proof}

\section{Width}\label{width}

In this section we prove the following result; its statement parallels that of Theorem
\ref{equality:depth}. It extends the result of Frankild \cite{Fr} that treats the
case of complexes with bounded below homology.

\begin{theorem}\label{equality:width}
Let $R$ be a noetherian ring and $M$ a complex of $R$-modules. Let ${\fa}$ be an ideal
in $R$ and $K$ the Koszul complex on a finite generating sequence for ${\fa}$. Then
\[
\inf(K\otimes_RM) = \inf(\Ltor{(R/{\fa})}RM) = \inf\Lho {\fa}M\,.
\]
\end{theorem}

The proof becomes more transparent with the following result on hand.

\begin{proposition}\label{weak sensitivity:width}
Let $M$ and $X$ be complexes of $R$-modules. If\/ $\hh X$ is bounded and ${\fa}^d\5\hh
X=0$ for some positive integer $d$, then
\[
\inf\left(\Ltor XRM\right)\geq \inf\left(\Ltor{R/{\fa}}RM\right) + \inf X\,.
 \]
\end{proposition}
\begin{proof}
Let $E$ be a faithfully injective $R$-module. Then from Proposition 
\ref{weak sensitivity:depth}, applied to the complexes of $R$-modules $\Rhom RME$ 
and $X$, and adjointness one obtains that
\[
\sup\Rhom R{\Ltor XRM}E \leq \sup\Rhom R{\Ltor{R/\fa}RM}E -\inf X\,.
 \]
This gives the desired result as $\sup\Rhom RLE=-\inf L$ for any complex $L$. 
\end{proof}

\begin{proof}[Proof of Theorem \ref{equality:width}]
As before, we identify $\Ltor KRM$ and $K\otimes_RM$.

 {\em Proof of\/} $\inf\left(\Ltor KRM\right) = \inf\left(\Ltor
{(R/{\fa})}RM\right)$. Let $E$ be a faithfully injective $R$-module. Then $\inf L = -
\sup\Rhom RLE$ for any complex $L$; thus, by adjointness, the desired equality is
equivalent to
\[
-\sup\Rhom RK{\Rhom RME} = -\sup\Rhom R{R/{\fa}}{\Rhom RME}\,.
\]
This is a special case of \eqref{equality:depth}.

{\em Proof of \/} $\inf(K\otimes_RM) = \inf\Lho {\fa}M$. We may assume that one of the
two quantities in consideration is not $-\infty$. Suppose that $\inf\Lho {\fa}M
>-\infty$. By \eqref{identity:lambda}, since ${\fa}\5\hh{K\otimes_RM}=0$, one has the the
first of the following isomorphisms, whilst the second is a special case of
\eqref{commutes:lambda}, since $K$ is finite and free.
\[
K\otimes_RM \simeq \Lho {\fa}{K\otimes_RM} \simeq K\otimes_R \Lho {\fa}M\,.
\]
This entails $\inf(K\otimes_RM)=\inf(\Lho {\fa}M)$, by \eqref{NAK:lambda} and
(\ref{accounts}.1).

Suppose now that $\inf\left(K\otimes_RM\right)=i>-\infty$. We know now that this is
equivalent to $\inf\left(\Ltor{R/{\fa}}RM\right)=i$. By \eqref{weak sensitivity:width},
this implies that
\begin{equation*}
\inf\left(\Ltor{R/{\fa}^n}RM\right)\geq i\qquad\text{for each integer $n\geq
1$}\,.\tag{$\dagger$}
  \end{equation*}
By definition, $\Lho {\fa}M \simeq \Lambda_{\fa}(F)$, where $F$ is an appropriate flat
resolution of $M$. The complex $\Lambda_{\fa}(F)$ is defined by the exactness of the
sequence of complexes
\[
0\to\Lambda_{\fa}(F)\to \prod ((R/{\fa}^n)\otimes_RF)
   \arto{\theta}\prod((R/{\fa}^n)\otimes_RF) \dashrightarrow 0\,,
\]
where the degree $\ell\,$th component of $\prod((R/{\fa}^n)\otimes_RF)$ is
$\prod((R/{\fa}^n)\otimes_RF_\ell)$, and $\theta(f_n)=(f_n - \pi^{n+1}(f_{n+1}))$, with
$\pi^{n+1}\col(R/{\fa}^{n+1})\otimes_RF \to (R/{\fa}^n)\otimes_RF$ the canonical
surjection.  The inverse system $\{\,(R/{\fa}^n)\otimes_RF\,\}_{n\geq 0}$ being
surjective, the map $\theta$ in the diagram above is surjective. For each integer $n$,
since $(R/{\fa}^n)\otimes_RF\simeq \Ltor {R/{\fa}^n}RM$, inequality ($\dagger$) ensures
that $\inf((R/{\fa}^n)\otimes_RF)$ is at least $i$. Therefore, from the homology long
exact sequence arising from the short exact sequence of complexes above we deduce that
$\inf\Lambda_{\fa}(F)\geq i-1$, that is to say, $\inf\Lho {\fa}M\geq i-1$. In particular,
$\inf\Lho {\fa}M>-\infty$, which brings us back to preceding case.

 This completes the proof of the theorem.
\end{proof}

\begin{definition}\label{definition:width}
Let ${\fa}$ be an ideal in $R$ and $K$ the Koszul complex on ${\fa}$. For a complex of
$R$-modules $M$, the {\em ${\fa}$-width of $M$ over $R$} is defined by one of the
following equivalent formulas:
\[
\Rwidth R{\fa}M = \begin{cases}
   \inf(K\otimes_RM)\,; \cr
   \inf\{\,\ell\in\BZ\mid \Tor \ell R{R/{\fa}}M\ne 0\,\}\,; \cr
   \inf\{\,\ell\in\BZ\mid \lho \ell {\fa}M\ne 0\,\}\,.
\end{cases}
\]
For a complex $M$ over a local ring $(R,\fm,k)$, the {\em width of $M$} is the number
\[
\mwidth_R M =\Rwidth R{\fm}M\,.
\]
 \end{definition}

For complexes whose homology is bounded below, width was introduced by Yassemi in
\cite{Ya} over local rings. In full generality, width for complexes was studied by
Christensen, Foxby and Frankild, who defined it via Kozsul homology \cite{CFF}.

Next we establish a tight relation between ${\fa}$-depth and ${\fa}$-width.

\begin{proposition}\label{depth-width}
Let ${\fa}$ be an ideal in $R$ and $M$ a complex of $R$-modules. Let $E$ a faithful
injective $R$-module and set $M^\vee=\Rhom RME$. Then
\[
\Rdepth R{\fa}M = \Rwidth R{\fa}{M^\vee}\qnd \Rwidth R{\fa}M = \Rdepth
R{\fa}{M^\vee}\,.
\]
\end{proposition}
\begin{proof}
Let $K$ be the Koszul complex on a finite sequence of generators for ${\fa}$. Since $K$ is
bounded and degreewise finite and free, $K\otimes_R{M^\vee}$ and ${\Hom RKM}^\vee$ are
isomorphic (see, for example, the proof of \cite[(4.4)]{AF}\,). This gives the second
equality below, whereas the first is due to the equality $\sup L=-\inf(L^\vee)$, which
holds for any complex $L$.
\[
-\sup\Hom RKM = \inf\Hom RKM^\vee = \inf\left(K\otimes_R{M^\vee}\right)\,.
\]
So $\Rdepth R{\fa}M = \Rwidth R{\fa}{M^\vee}$; the proof of the other equality is
analogous.
\end{proof}

The preceding proposition allows us to obtain many a results concerning width from
those on depth. For instance, the following observation is a immediate consequence of
\eqref{infinity:depth}.

\begin{noname}\label{infinity:width}
Let ${\fa}$ be an ideal in $R$ and let $M$ be a complex of $R$-modules with $\hh M$
degreewise finite. Then
\begin{alignat*}{2}
&\Rwidth R{\fa}M = \infty& &\iff \inf M = \infty\,;\\
&\Rwidth R{\fa}M = -\infty & &\iff\inf M = -\infty\,.
 \end{alignat*}
\end{noname}

One can translate also \eqref{AB:depth} - the Auslander-Buchsbaum formula - into a
statement concerning width. We end this section with the formula for the width of the
derived tensor product, and the derived homomorphisms, of a pair of complexes.

\begin{proposition}\label{AB:width}
Let $(R,\fm,k)$ be a noetherian local ring, and let $M$ and $N$ be complexes of
$R$-modules.
\begin{align*}
\mwidth_R(\Ltor MRN) &= \mwidth_RM + \mwidth_RN\,;\\
\mdepth_R(\Rhom RMN) & = \mwidth_RM + \mdepth_RN \,.
\end{align*}
\end{proposition}
\begin{proof}
The associativity formula $\Ltor{(\Ltor MRN)}Rk \simeq (\Ltor MRk)\otimes_k(\Ltor kRN)$
leads to the first equality below, while the second holds since $k$ is a field.
\[
\inf\left(\Ltor{(\Ltor MRN)}Rk\right)= \inf\left((\Ltor MRk)\otimes_k(\Ltor kRN)\right)
          = \inf(\Ltor MRk) + \inf(\Ltor kRN)\,.
\]
This establishes the first formula; the argument for the second one is similar.
\end{proof}

\end{document}